\documentclass[11pt]{article}
\usepackage{amsfonts}
\usepackage{latexsym, amssymb, amsmath, amscd, amsfonts, epsfig, graphicx, colordvi}
\usepackage{ifpdf}
\usepackage{graphicx}
\usepackage{xcolor}
\usepackage[dvipsnames]{pstricks}
\newtheorem{thm}{Theorem}[section]

\newtheorem{lem}[thm]{Lemma}

\def\pf{\noindent{\it Proof.} }
\def\qed{\nopagebreak\hfill{\rule{4pt}{7pt}}
\medbreak}

\def\qed{\nopagebreak\hfill{\rule{4pt}{7pt}}
\medbreak}

\makeatletter
\def\ExtendSymbol#1#2#3#4#5{\ext@arrow 0099{\arrowfill@#1#2#3}{#4}{#5}}
\makeatother

\title{Congruences modulo $4$ for Rogers--Ramanujan--Gordon type overpartitions}
\author{
Doris D. M. Sang\raisebox{5pt}{\scriptsize 1} and Diane Y. H.
Shi\raisebox{5pt}{\scriptsize 2}}
\date{
\vspace{15pt}\raisebox{5pt}{\scriptsize 1\,}School of Mathematics
and Quantitative Economics\\Dongbei
University of Finance and Economics, Liaoning 116025, P.R. China\\sdm@cfc.nankai.edu.cn\\
\vspace{15pt}\raisebox{5pt}{\scriptsize 2\,}School of
Mathematics\\ Tianjin University, Tianjin 300072, P.R.
China\\shiyahui@tju.edu.cn}

\begin{document}
\maketitle

\vspace{0.3cm} \noindent{\bf Abstract.}
In a  recent work, Andrews defined the singular overpartitions with the goal of presenting an overpartition analogue to the theorems of Rogers--Ramanujan type for ordinary partitions with restricted successive ranks. As a small part of his work, Andrews noted two congruences modulo $3$ for the number of singular overpartitions prescribed by parameters $k=3$ and $i=1$. It should be noticed that this number equals the number of the Rogers--Ramanujan--Gordon type overpartitions with $k=i=3$ which come from the overpartition analogue of Gordon's Rogers--Ramanujan partition theorem introduced by Chen, Sang and Shi. In this paper, we  derive numbers of  congruence identities modulo $4$ for the number of Rogers--Ramanujan--Gordon type overpartitions.



\section{Introduction}
In a  recent work, Andrews \cite{and15} introduced  the singular overpartitions to present a Rogers--Ramanujan type theorem for overpartitions and denoted $\overline{Q}_{k,i}(n)$ to be the number of singular overpartitions of $n$ subject to overlining conditions prescribed by $k$ and $i$.
As a part of his work, Andrews noted two congruences modulo $3$  for $\overline{Q}_{k,i}(n)$ with $k=3$ and $i=1$, which equals  the number of overpartitions of $n$ into parts not  divisible by $3$. This number  is also equal to   the number $\overline{A}_{3,3,1}(n)$, where $\overline{A}_{k,i,1}(n)$ is defined to be the number of Rogers--Ramanujan--Gordon type overpartitions of $n$ with even moduli.

The aim of this paper is to derive  some congruences for the number of Rogers--Ramanujan-Gordon type overpartitions
modulo $4$. We shall prove a number of results by constructing 3-dissection and $4$-dissection of the generating function for
the Rogers--Ramanujan--Gordon type overpartitions.

Let us give an overview of  some definitions. A partition $\lambda$
of a positive integer $n$ is a non-increasing sequence of positive
integers $\lambda_1\geq \cdots\geq \lambda_s>0$ such that
$n=\lambda_1+\cdots+\lambda_s$. The partition of zero is the
partition with no parts. An overpartition $\lambda$ of a positive
integer $n$ is also a non-increasing sequence of positive integers
$\lambda_1\geq \cdots\geq \lambda_s>0$ such that
$n=\lambda_1+\cdots+\lambda_s$ and the first occurrence  of  each
integer may be overlined. We denote the number of overpartitions
of $n$ by $\overline{p}(n)$. For example, there are $8$ overpartitions of $3$:
\[3;\ \overline{3};\ 2+ 1;\ \overline{2} + 1;\ 2 + \overline{1};\ \overline{2} + \overline{1};\ 1 + 1 + 1;\ \overline{1} + 1 + 1.\]
For a partition
or an overpartition $\lambda$ and for any integer $l$, let
$f_l(\lambda) (f_{\overline{l}}(\lambda))$ denote the number of
occurrences of non-overlined (resp.~overlined) $l$ in $\lambda$. Let
$V_{\lambda}(l)$ denote the number of overlined parts in $\lambda$
that are less than or equal to $l$.
We shall adopt  the common notation  as used in Andrews \cite{and76}. Let
 \[(a)_\infty=(a;q)_\infty=\prod_{i=0}^{\infty}(1-aq^i),\]
 and \[(a)_n=(a;q)_n=\frac{(a)_\infty}{(aq^n)_\infty}.\]
We also write
 \[(a_1,\ldots,a_k;q)_\infty=(a_1;q)_\infty\cdots(a_k;q)_\infty.\]

 The Rogers-Ramanujan-Gordon theorems for overpartition are motivated by the Rogers--Ramanujan--Gordon type partition. The first combinatorial generalization of the Rogers--Ramanujan identities was derived by Gordon in 1961 \cite{gor61}. And then, in 1967,  Andrews \cite{and66} gave the analytic generalization of Rogers--Ramanujan identities, which is also a generating function form of Gordon's theorem.  These two theorems  both concern odd moduli. In 1980£¬ Bressoud \cite{Bre79} extended their results to all moduli.

\begin{thm}\label{Bressoud}
Given positive integer $k$, $j = 0$ or $1$ and integral $i$ such that $0<i<(2k+j)/2$. Let $A_{k,i,j}(n)$ denote the number of partitions of $n$ that no part is congruent to $0$ or $\pm i$ modulo $2k+j$. Let  $B_{k,i,j}(n)$ denote the number of
partitions of $n$ of the form
$\lambda_1+\lambda_2+\cdots+\lambda_s$ such that:
\begin{itemize}
\item[(i)]$f_1(\lambda)\leq i-1$;
\item[(ii)] $f_l(\lambda)+f_{l+1}(\lambda)\leq k-1$;
\item[(iii)] if $f_l(\lambda)+f_{l+1}(\lambda)=k-1$, then
$lf_l(\lambda)+(l+1)f_{l+1}(\lambda)  \equiv i-1 \pmod{2-j}$.
\end{itemize}
Then we have
\begin{equation}
A_{k,i,j}(n)=B_{k,i,j}(n).
\end{equation}
\end{thm}

Bressoud \cite{Bre80} also gave the analytic form of this theorem.
 \begin{thm} Given positive integer $k$, $j = 0$ or $1$ and integral $i$ such that $0<i<(2k+j)/2$,
\begin{equation}\label{AB}\sum_{N_1\geq N_2\geq\cdots\geq
N_{k-1}\geq0}\frac{q^{N_1^2+N_2^2+\cdots+N_{k-1}^2+N_{i+1}+\cdots+N_{k-1}}
}{(q)_{N_1-N_2}\cdots(q)_{N_{k-2}-N_{k-1}}(q^{2-j};q^{2-j})_{N_{k-1}}}=\frac{(q^i,q^{2k+j-i},q^{2k+j};q^{2k+j})_\infty}{(q)_\infty}.
\end{equation}
\end{thm}

One can see that, when $j=1$ Theorem \ref{Bressoud} becomes the Gordon's generalization in \cite{gor61} and identity \eqref{AB} becomes the analytic form of
Andrews obtained in \cite{and66}.

In 2013,  Chen, Sang and Shi \cite{chen13} found  an overpartition analogue of the Rogers--Ramanujan--Gordon theorem,
 more precisely, an overpartition analogue of Theorem \ref{Bressoud} in the case $j=1$ and also the generating function form.
 In 2015, Chen, Sang and Shi \cite{chen15} obtained the overpartition analogue of Theorem \ref{Bressoud} in the case $j=0$. The generating function form was also derived by Sang and Shi \cite{sang15}. Here we state these two theorems in a unified form  together with its generating function form.

\begin{thm}\label{main1}Let $j=0$ or $1$, and $k$ and $i$ be integers, such that $0<i<(2k+j)/2$. Let $\overline{B}_{k,i,j}(n)$ denote the number of overpartitions of $n$ of the form $\lambda_1+\lambda_2+\cdots+\lambda_s$,
such that
\begin{itemize}\item[(i)] $f_1(\lambda)\leq i-1$;
\item[(ii)] $f_l(\lambda)+f_{\overline{l}}(\lambda)+f_{l+1}(\lambda)\leq k-1$;
\item[(iii)] if the equality in Condition (ii) is attained at $l$,
 i.e., $f_l(\lambda)+f_{\overline{l}}(\lambda)+f_{l+1}(\lambda)=k-1$, then $lf_l(\lambda)+lf_{\overline{l}}(\lambda)+(l+1)f_{l+1}(\lambda)\equiv V_{\lambda}(l)+i-1 \pmod{2-j}$.\end{itemize}
 For $i\neq k$ nor $j\neq1$, let $\overline{A}_{k,i,j}(n)$ denote the number of
overpartitions of $n$ whose non-overlined parts are not congruent to
$0,\pm i$ modulo $2k-1+j$. Let $\overline{A}_{k,k,1}(n)$ denote the number of overpartitions of $n$ with parts not divisible by $k$. Then for all $n\geq 0$, we have
\begin{equation}\overline{A}_{k,i,j}(n)=\overline{B}_{k,i,j}(n).
\end{equation}
 \end{thm}

The generating function form of
  Theorem \ref{main1}, which is an
Andrews--Gordon type identity for overpartitions with all moduli can be stated as follows.

 \begin{thm}\label{CDthm}
Let $j=0$ or $1$ and  $k$ and $i$ be integers, such that $0<i<(2k+j)/2$, we have
\begin{align}\label{CD}
&\sum_{N_1\geq N_2\geq\cdots\geq N_{k-1}\geq0}
\frac{q^{\frac{(N_1+1)N_1}{2}+N_2^2+\cdots+N_{k-1}^2+N_{i+1}+\cdots+N_{k-1}}
(-q)_{N_1-1}(1+q^{N_i})}{(q)_{N_1-N_2}\cdots(q)_{N_{k-2}-N_{k-1}}
(q^{2-j};q^{2-j})_{N_{k-1}}}\nonumber
\\[6pt]&\qquad =\frac{(-q)_\infty(q^i,q^{2k-1-i+j},q^{2k-1+j};q^{2k-1+j})_\infty}{(q)_\infty}.
\end{align}
\end{thm}

In this paper, we shall derive a number of congruence properties of $A_{k,i,j}(n)$, which can be called the number of overpartitions of Rogers--Ramanujan--Gordon type. Indeed, two congruence properties for a number equal to $A_{3,3,1}(n)$ have been obtained by Andrews.

In 2015, Andrews \cite{and15} defined a new type of  overpartitions,  called singular overpartition to give the overpartition analogous to Rogers--Ramanujan type theorems for ordinary
partitions with restricted successive ranks and  denoted $\overline{Q}_{k,i}(n)$ to be the number of singular overpartitions of $n$ subject
to overlining conditions prescribed by $k$ and $i$.  The generating function of $\overline{Q}_{k,i}(n)$ was derived as follows:
\begin{equation}
\sum_{n=0}^\infty \overline{Q}_{k,i}(n)q^n=\frac{(q^k,-q^i,-q^{k-i};q^k)_\infty}{(q;q)_\infty}.
\end{equation}
As a small part of his work, Andrews gave the following result:
\begin{thm}
\begin{equation}\overline{Q}_{3,1}(9n+3)\equiv\overline{Q}_{3,1}(9n+3)\equiv0\pmod{3}.
\end{equation}
\end{thm}

It  should be noticed that
\[\overline{Q}_{3,1}(n)=\overline{B}_{3,3,1}(n)=\overline{A}_{3,3,1}(n), \]
that is to say,
\begin{equation}\nonumber
\overline{A}_{3,3,1}(9n+3)\equiv\overline{A}_{3,3,1}(9n+3)\equiv0\pmod{3}.
\end{equation}

We shall consider more congruence properties of $\overline{A}_{k,i,j}(n)$. To make the computation easier, we define $S_{k,i}(n)$ as follows.
\begin{equation}\nonumber S_{2k-1+j,i}(n)=\overline{A}_{k,i,j}(n),
\end{equation}
then for $1\leq i\leq k/2$,
\begin{equation}\label{s}\sum_{n\geq 0}S_{k,i}(n)q^n=\frac{(-q)_\infty(q^i,q^{k-i},q^{k};q^{k})_\infty}{(q)_\infty}.
\end{equation}

The proofs of the congruences of $S_{k,i}(n)$  involve the following definitions and results.
Let us employ the Ramaujan's general theta function $f(a,b)$  defined by
\begin{equation}f(a,b):=\sum_{n=-\infty}^{\infty}a^{\frac{n(n+1)}{2}}b^{\frac{n(n-1)}{2}},\qquad |ab|<1.
\end{equation}
One special case is, in Ramanujan's notation,
\begin{equation}\label{eq:DEFf}
f(-q)=f(-q,-q^2)=\sum_{n=-\infty}^{\infty}(-1)^nq^{n(3n-1)/2}.
\end{equation}
By using Jacobi's triple product identity that
\begin{equation}\label{J3}
\sum_{n=-\infty}^{\infty}z^nq^{n^2}=(-zq;q^2)_\infty(-q/z;q^2)_\infty(q^2;q^2)_\infty,
\end{equation}
it follows
\begin{equation}\nonumber
f(-q)=(q;q)_\infty.
\end{equation}
Throughout this paper, we use $f_n$ to denote $f(-q^n)$, that is,
\begin{equation}
f_n:=f(-q^n)=(q^n;q^n)_\infty.
\end{equation}

As noted by Corteel and Lovejoy \cite{cor04}, the generating function of $\overline{p}(n)$ is given by
\begin{equation}
\sum_{n=0}^\infty \overline{p}(n)q^n=\frac{(-q;q)_\infty}{(q)_\infty}=\frac{f_2}{f_1^2}.
\end{equation}

It can be  seen that the generating function $S_{k,i}(n)$ can be written in terms of Ramanujan's  theta function as follows.
 \begin{equation}\sum_{n\geq 0}S_{k,i}(n)q^n=\frac{(-q)_\infty(q^i,q^{k-i},q^{k};q^{k})_\infty}{(q)_\infty}
 =f(-q^i,-q^{k-i})\frac{f_2}{f_1^2}.
\end{equation}

In \cite{Hir05}, Hirschhorn and Sellers  obtained $2$-, $3$- and $4$-dissections of the
generating function of $\overline{p}(n)$ and derived a number of congruences for $\overline{p}(n)$ modulo $4$, $8$
and $64$ including $\overline{p}(8n+7) \equiv0 \pmod{64}$, for $n \geq 0$.

In this paper we shall use  the following $3$-dissection and $4$-dissection of the generating function of $\overline{p}(n)$ given by Hirschhorn and Seller \cite{Hir05} to compute the congruences of $S_{k,i}(n)$.
\begin{thm}
The $3$-dissection of the generating function of $\overline{p}(n)$ is
\begin{equation}\label{p3}
\sum_{n=0}^\infty \overline{p}(n)q^n
=\frac{f_2}{f_1^2}=\frac{f_6^4f_9^6}{f_3^8f_{18}^3}
+2q\frac{f_6^3f_9^3}{f_3^7}
+4q^2\frac{f_6^2f_{18}^3}{f_3^6},
\end{equation}
and the $4$-dissection is

\begin{equation}\label{p4}
\sum_{n=0}^\infty \overline{p}(n)q^n
=\frac{f_8^{19}}{f_4^{14}f_{16}^6}
+2q\frac{f_8^{13}}{f_4^{12}f_{16}^2}
+4q^2\frac{f_8^7f_{16}^2}{f_4^{10}}
+8q^3\frac{f_8f_{16}^6}{f_4^8}.
\end{equation}
\end{thm}

Our main goal in this paper is to prove numerous arithmetic relations
satisfied by $S_{k,i}(n)$.  The techniques we employ are
elementary, by involving dissections of $q$-series. In Section 2, we shall compute the $3$-dissection and $4$-dissection of $f(-q^i,-q^{k-i})$.
In Section 3 and 4, we will discuss the congruences of $S_{k,i}(n)$ by considering $3$-dissection and $4$-dissection of the generating function of $S_{k,i}(n)$.

\section{The dissections of $f(-q^i,-q^{k-i})$}

In this section, by using Jacobi's triple product identity \eqref{J3} we give the $3$-dissection and $4$-dissection of $f(-q^i,-q^{k-i})$.
\begin{lem}
The $3$-dissection of $f(-q^i,-q^{k-i})$ is as follows
\begin{equation}f(-q^i,-q^{k-i})
\label{d3}=f(-q^{3i+3k},-q^{6k-3i})-q^if(-q^{3i+6k},-q^{3k-3i})+q^{2i+k}f(-q^{-3i},-q^{3i+9k}).
\end{equation}
\end{lem}
\pf By definition \eqref{eq:DEFf},
 \[f(-q^i,-q^{k-i})
=\sum_{n=-\infty}^{\infty}(-1)^nq^{in}q^{kn(n-1)/2}=\sum_{n=-\infty}^{\infty}(-1)^nq^{[kn^2+(2i-k)n]/2},\]
so that  the $3$-dissection is
\begin{align*}f(-q^i,-q^{k-i})
&=\sum_{n=-\infty}^{\infty}(-1)^{3n}q^{[9kn^2+3(2i-k)n]/2}
  +\sum_{n=-\infty}^{\infty}(-1)^{3n+1}q^{[9kn^2+(6i+3k)n+2i]/2}
  \\&=\sum_{n=-\infty}^{\infty}(-1)^{3n+2}q^{[9kn^2+(6i+9k)n+(4i+2k)]/2}.
\end{align*}
  By employing Jacobi's triple identity,  we have
  \begin{align*}
f(-q^i,-q^{k-i})=&(q^{3i+3k},q^{6k-3i},q^{9k};q^{9k})_\infty-q^i(q^{3i+6k},q^{3k-3i},q^{9k};q^{9k})_\infty\\
&+q^{2i+k}(q^{-3i},q^{3i+9k},q^{9k};q^{9k})_\infty
\\=&f(-q^{3i+3k},-q^{6k-3i})-q^if(-q^{3i+6k},-q^{3k-3i})+q^{2i+k}f(-q^{-3i},-q^{3i+9k}).
\end{align*}\qed

 Similarly as the $3$-dissection of $f(-q^i,-q^{k-i})$,   the following $4$-dissection of $f(-q^i,-q^{k-i})$  follows.
\begin{lem} We have
\begin{align}\nonumber
f(-q^i, -q^{k-i})=&f(-q^{6k+4i},-q^{10k-4i})
-q^if(-q^{10k+4i},-q^{6k-4i})
\\&\label{diss4}+q^{2i+k}f(-q^{14k+4i},-q^{2k-4i})
-q^{3i+3k}f(-q^{18k+4i},-q^{-2k-4i}).
\end{align}
\end{lem}
\pf
\begin{align*}
f(-q^i, -q^{k-i})
=&\sum_{n=-\infty}^{\infty}(-1)^nq^{[kn^2+(2i-k)n]/2}
\\=&\sum_{n=-\infty}^{\infty}q^{8kn^2+2(2i-k)n}
-q^i\sum_{n=-\infty}^{\infty}q^{8kn^2+(4i+2k)n}
\\&+q^{2i+k}\sum_{n=-\infty}^{\infty}q^{8kn^2+(6k+4i)n}
-q^{3i+3k}\sum_{n=-\infty}^{\infty}q^{8kn^2+(10k+4i)n}
\\=&f(-q^{6k+4i},-q^{10k-4i})
-q^if(-q^{10k+4i},-q^{6k-4i})
\\&+q^{2i+k}f(-q^{14k+4i},-q^{2k-4i})
-q^{3i+3k}f(-q^{18k+4i},-q^{-2k-4i}).
\end{align*}\qed

In next section, we shall consider this $3$-dissection of the generating function of $S_{k,i}(n)$  corresponding to the  parameters $k$ and $i${\color{red},} and get  arithmetic properties of $S_{k,i}(n)$ by employing the $3$-dissection of $\overline{p}(n)$.

\section{Arithmetic properties as the consequence of the $3$-dissection of the generating function of $S_{k,i}(n)$}

Combining \eqref{d3} and \eqref{p3}, we get the following $3$-dissection of the generating function of $S_{k,i}(n)$:

\begin{align}\nonumber\sum_{n\geq 0}S_{k,i}(n)q^n=&f(-q^i,-q^{k-i})\frac{f_2}{f_1^2}
\\ \nonumber=&[f(-q^{3i+3k},-q^{6k-3i})-q^if(-q^{3i+6k},-q^{3k-3i})+q^{2i+k}f(-q^{-3i},-q^{3i+9k})]
\\\nonumber&\times\left(\frac{f_6^4f_9^6}{f_3^8f_{18}^3}+2q\frac{f_6^3f_9^3}{f_3^7}
+4q^2\frac{f_6^2f_{18}^3}{f_3^6}\right)
\\ \nonumber=&f(-q^{3i+3k},-q^{6k-3i})\frac{f_6^4f_9^6}{f_3^8f_{18}^3}
-q^if(-q^{3i+6k},-q^{3k-3i})\frac{f_6^4f_9^6}{f_3^8f_{18}^3}
+q^{2i+k}f(-q^{-3i},-q^{3i+9k})\frac{f_6^4f_9^6}{f_3^8f_{18}^3}
\\&\nonumber+2qf(-q^{3i+3k},-q^{6k-3i})\frac{f_6^3f_9^3}{f_3^7}
-2q^{i+1}f(-q^{3i+6k},-q^{3k-3i})\frac{f_6^3f_9^3}{f_3^7}
\\ \nonumber&+2q^{2i+k+1}f(-q^{-3i},-q^{3i+9k})\frac{f_6^3f_9^3}{f_3^7}+4q^2f(-q^{3i+3k},-q^{6k-3i})\frac{f_6^2f_{18}^3}{f_3^6}
\\&\label{S3}
-4q^{i+2}f(-q^{3i+6k},-q^{3k-3i})\frac{f_6^2f_{18}^3}{f_3^6}
+4q^{2i+k+2}f(-q^{-3i},-q^{3i+9k})\frac{f_6^2f_{18}^3}{f_3^6}.
\end{align}

\subsection{The congruences of $S_{k,i}(n)$ in the case of  $k\equiv2\pmod{3}$, $i\equiv1\pmod{3}$}

In the case of $k\equiv2\pmod{3}$ and  $i\equiv1\pmod{3}$ we rearrange \eqref{S3} to get the following $3$-dissection :
\begin{align*}&\sum_{n\geq 0}S_{k,i}(n)q^n=f(-q^i,-q^{k-i})\frac{f_2}{f_1^2}
\\=&\left[f(-q^{3i+3k},-q^{6k-3i})\frac{f_6^4f_9^6}{f_3^8f_{18}^3}
-4q^{i+2}f(-q^{3i+6k},-q^{3k-3i})\frac{f_6^2f_{18}^3}{f_3^6}
+4q^{2i+k+2}f(-q^{-3i},-q^{3i+9k})\frac{f_6^2f_{18}^3}{f_3^6}\right]
\\&+q\left[-q^{i-1}f(-q^{3i+6k},-q^{3k-3i})\frac{f_6^4f_9^6}{f_3^8f_{18}^3}
+q^{2i+k-1}f(-q^{-3i},-q^{3i+9k})\frac{f_6^4f_9^6}{f_3^8f_{18}^3}
+2f(-q^{3i+3k},-q^{6k-3i})\frac{f_6^3f_9^3}{f_3^7}\right]
\\&+q^2\left[-2q^{i-1}f(-q^{3i+6k},-q^{3k-3i})\frac{f_6^3f_9^3}{f_3^7}
+2q^{2i+k-1}f(-q^{-3i},-q^{3i+9k})\frac{f_6^3f_9^3}{f_3^7}
+4f(-q^{3i+3k},-q^{6k-3i})\frac{f_6^2f_{18}^3}{f_3^6}\right].
\end{align*}
On the right-hand side of the above identity in each  square brackets, the powers of $q$ are all multiples of $3$.

Then, we can get the generating function of $S_{k,i}(3n)$ as follows

\begin{equation}\nonumber
\sum_{n\geq 0}S_{k,i}(3n)q^n
=f(-q^{i+k},-q^{2k-i})\frac{f_2^4f_3^6}{f_1^8f_{6}^3}
-4q^{\frac{i+2}{3}}f(-q^{i+2k},-q^{k-i})\frac{f_2^2f_{6}^3}{f_1^6}+4q^{\frac{2i+k+2}{3}}f(-q^{-i},-q^{i+3k})\frac{f_2^2f_{6}^3}{f_1^6},
\end{equation}
which implies that
\begin{equation}\nonumber
\sum_{n\geq 0}S_{k,i}(3n)q^n\equiv f(-q^{i+k},-q^{2k-i})\frac{f_2^4f_3^6}{f_1^8f_{6}^3}\pmod{4}.
\end{equation}
Applying \eqref{p3} to $f_2^4/f_1^8$, we have
\begin{equation}\label{q24q8}
\frac{f_2^4}{f_1^8}=\left(\frac{f_2}{f_1^2}\right)^4
=\left(\frac{f_6^4f_9^6}{f_3^8f_{18}^3}+2q\frac{f_6^3f_9^3}{f_3^7}
+4q^2\frac{f_6^2f_{18}^3}{f_3^6}\right)^4
\equiv\left(\frac{f_6^4f_9^6}{f_3^8f_{18}^3}\right)^4\pmod{4}.
\end{equation}

So we have
\begin{equation}\label{kis213}\sum_{n\geq 0}S_{k,i}(3n)q^n\equiv f(-q^{i+k},-q^{2k-i})\left(\frac{f_6^4f_9^6}{f_3^8f_{18}^3}\right)^4\frac{f_3^6}{f_{6}^3}\pmod{4}.
\end{equation}
It can be checked  that all exponents of $q$ in the right-hand side of \eqref{kis213} are multiples of $3$. Then one can verify that
\begin{equation}\nonumber
S_{k,i}(9n+3)\equiv S_{k,i}(9n+6)\equiv 0(mod\ 4),
\end{equation}
with $k\equiv2\pmod{3}$, $i\equiv1\pmod{3}$.

\subsection{$k\equiv2\pmod{3}$, $i\equiv2\pmod{3}$}
In the case $k\equiv2\pmod{3}$ and  $i\equiv2\pmod{3}$ we rearrange \eqref{S3} to get the following $3$-dissection:

\begin{align*}
&\sum_{n\geq 0}S_{k,i}(n)q^n=f(-q^i,-q^{k-i})\frac{f_2}{f_1^2}
\\=&\left[f(-q^{3i+3k},-q^{6k-3i})\frac{f_6^4f_9^6}{f_3^8f_{18}^3}
+q^{2i+k}f(-q^{-3i},-q^{3i+9k})\frac{f_6^4f_9^6}{f_3^8f_{18}^3}
-2q^{i+1}f(-q^{3i+6k},-q^{3k-3i})\frac{f_6^3f_9^3}{f_3^7}\right]
\\&+q\left[2f(-q^{3i+3k},-q^{6k-3i})\frac{f_6^3f_9^3}{f_3^7}
+2q^{2i+k}f(-q^{-3i},-q^{3i+9k})\frac{f_6^3f_9^3}{f_3^7}
-4q^{i+1}f(-q^{3i+6k},-q^{3k-3i})\frac{f_6^2f_{18}^3}{f_3^6}\right]
\\&+q^2\left[-q^{i-2}f(-q^{3i+6k},-q^{3k-3i})\frac{f_6^4f_9^6}{f_3^8f_{18}^3}
+4f(-q^{3i+3k},-q^{6k-3i})\frac{f_6^2f_{18}^3}{f_3^6}
+4q^{2i+k}f(-q^{-3i},-q^{3i+9k})\frac{f_6^2f_{18}^3}{f_3^6}\right].
\end{align*}

By this $3$-dissection for $k\equiv2\pmod{3}$, $i\equiv2\pmod{3}$, we can get the generating function of $S_{k,i}(3n+2)$.

\begin{align}&\nonumber\sum_{n\geq 0}S_{k,i}(3n+2)q^n
\\=&4f(-q^{i+k},-q^{2k-i})\frac{f_2^2f_6^3}{f_1^6}
+4q^{(2i+k)/3}f(-q^{-i},-q^{i+3k})\frac{f_2^2f_6^3}{f_1^6}
-q^{(i-2)/3}f(-q^{i+2k},-q^{k-i})\frac{f_2^4f_3^6}{f_1^8f_6^3},
\end{align}
and  the following arithmetic property, that,
\begin{equation}
\sum_{n\geq 0}S_{k,i}(3n+2)q^n\equiv-q^{(i-2)/3}f(-q^{i+2k},-q^{k-i})\frac{f_2^4f_3^6}{f_1^8f_6^3}\pmod{4}
\end{equation}
As a consequence of \eqref{q24q8}, we have
\begin{equation}
\sum_{n\geq0}S_{k,i}(3n+2)q^n\equiv-q^{(i-2)/3}f(-q^{i+2k},-q^{k-i})
\left(\frac{f_6^4f_9^6}{f_3^8f_{18}^3}\right)^4\frac{f_3^6}{f_6^3}\pmod{4}.
\end{equation}
Then, for $k\equiv 2\pmod{3}$, we have
\begin{itemize}
\item[1.]if $i\equiv 2\pmod{9}$, then \begin{equation}S_{k,i}(9n+5)\equiv S_{k,i}(9n+8)\equiv0\pmod{4};\end{equation}
\item[2.]if $i\equiv 5\pmod{9}$, then\begin{equation}S_{k,i}(9n+2)\equiv S_{k,i}(9n+8)\equiv0\pmod{4};\end{equation}
\item[3.]if $i\equiv 8\pmod{9}$, then \begin{equation}S_{k,i}(9n+2)\equiv S_{k,i}(9n+5)\equiv0\pmod{4}.\end{equation}
\end{itemize}

\subsection{$k\equiv0\pmod{3}$, $i\equiv0\pmod{3}$}
Similar with the above cases we can get the following identity
\begin{equation}\nonumber
\sum_{n=0}^\infty S_{k,i}(3n+2)q^n=4\frac{f_2^2f_6^3}{f_1^6}[f(-q^{i+k},-q^{2k-i})-q^{i/3}f(-q^{i+2k},-q^{k-i})
+q^{(2i+k)/3}f(-q^{-i},-q^{i+3k})].
\end{equation}
It is easy to see that, for $k\equiv i\equiv0\pmod{3}$, we have
\begin{equation}
S_{k,i}(3n+2)\equiv0\pmod{4}.
\end{equation}

\subsection{$k\equiv2\pmod{3}$, $i\equiv0\pmod{3}$}
In this case we also consider the generating function of $S_{k,i}(3n+2)$
\begin{equation}\nonumber
\sum_{n=0}^\infty S_{k,i}(3n+2)q^n=4[f(-q^{i+k},-q^{2k-i})-q^{i/3}f(-q^{i+2k},-q^{k-i})]
\frac{f_2^2f_3^3}{f_1^6}+q^{(2i+k-2)/3}f(-q^{-i},-q^{i+3k})\frac{f_2^4f_3^6}{f_1^8f_6^3}\end{equation}
So we have that
\begin{align}\nonumber
\sum_{n=0}^\infty S_{k,i}(3n+2)q^n
&\equiv q^{(2i+k-2)/3}f(-q^{-i},-q^{i+3k})\frac{f_2^4f_3^6}{f_1^8f_6^3}\pmod{4}\\\nonumber&\equiv q^{(2i+k-2)/3}f(-q^{-i},-q^{i+3k})\left(\frac{f_6^4f_9^6}{f_3^8f_{18}^3}\right)^4\frac{f_3^6}{f_6^3}\pmod{4}.
\end{align}
Then we can verify that
\begin{itemize}\nonumber
\item[1.]for $k-i\equiv 2\pmod{9}$,
\begin{equation}S_{k,i}(9n+5)\equiv S_{k,i}(9n+8)\equiv 0\pmod{4};
\end{equation}
\item[2.]for $k-i\equiv 5\pmod{9}$,
\begin{equation}S_{k,i}(9n+2)\equiv S_{k,i}(9n+8)\equiv 0\pmod{4};\end{equation}
\item[2.]for $k-i\equiv 8\pmod{9}$,
\begin{equation}S_{k,i}(9n+2)\equiv S_{k,i}(9n+5)\equiv 0\pmod{4}.\end{equation}
\end{itemize}

\section{Arithmetic properties as the consequence of the $4$-dissection of the generating function of $S_{k,i}(n)$}

In this section, we shall prove some arithmetic properties according to the $4$-dissection of the generating function of $S_{k,i}(n)$. The $4$-dissection of the generating function of $S_{k,i}(n)$ is based on the  $4$-dissection of $f(-q^i,-q^{k-i})$ given in \eqref{diss4}  and the $4$-dissection of $f_2/f_1^2$ given in \eqref{p4}.

Recall that  the $4$-dissection of $f_2/f_1^2$ is
\begin{equation}\nonumber
\frac{f_2}{f_1^2}
=\frac{f_8^{19}}{f_4^{14}f_{16}^6}
+2q\frac{f_8^{13}}{f_4^{12}f_{16}^2}
+4q^2\frac{f_8^7f_{16}^2}{f_4^{10}}
+8q^3\frac{f_8f_{16}^6}{f_4^8}.
\end{equation}

Combining \eqref{diss4} we have that
\begin{align}\nonumber
\nonumber \sum_{n=0}^{\infty}S_{k,i}(n)q^n\equiv &[f(-q^{6k+4i},-q^{10k-4i})-q^if(-q^{10k+4i},-q^{6k-4i})
\\&\nonumber+q^{2i+k}f(-q^{14k+4i},-q^{2k-4i})-q^{3i+3k}f(-q^{18k+4i},-q^{-2k-4i})]
\\&\times\left(\frac{f_8^{19}}{f_4^{14}f_{16}^6}+2q\frac{f_8^{13}}{f_4^{12}f_{16}^2}\right)\pmod{4}.
\end{align}
We can get the following results
\begin{itemize}
\item[1.]$k\equiv 2\pmod{4}$, $i\equiv1\pmod{4}$,
\begin{equation}
S_{k,i}(4n+3)\equiv0\pmod{4};
\end{equation}

\item[2.]$k\equiv 2\pmod{4}$, $i\equiv3\pmod{4}$,
\begin{equation}
S_{k,i}(4n+2)\equiv0\pmod{4};
\end{equation}
\item[3.]$k\equiv 3\pmod{4}$, $i\equiv0\pmod{4}$,
\begin{equation}
S_{k,i}(4n+3)\equiv0\pmod{4};
\end{equation}
\item[4.]$k\equiv 0\pmod{4}$, $i\equiv0\pmod{4}$,
\begin{equation}
S_{k,i}(4n+2)\equiv S_{k,i}(4n+3)\equiv0\pmod{4}.
\end{equation}
\end{itemize}

\vspace{0.5cm}
 \noindent{\bf Acknowledgments.}
This work was supported by the National Science Foundation of China (Nos.1140149, 11501089, 11501408).

\end{document}